\documentclass{amsart}

\usepackage{mathrsfs}

\usepackage[arrow,matrix]{xy}
\usepackage[all]{xypic}

\newtheorem{theorem}{Theorem}

\theoremstyle{definition}

\theoremstyle{remark}

\numberwithin{equation}{section}

\providecommand{\rad}{\mathop{\rm rad}\nolimits}
\renewcommand{\mod}{\mathop{\rm mod}\nolimits}%
\providecommand{\Tr}{\mathop{\rm Tr}\nolimits}%
\newcommand{\Hom}{\operatorname{Hom}}
\newcommand{\End}{\operatorname{End}}
\newcommand{\add}{\operatorname{add}}

\newcommand{\ann}{\operatorname{ann}}
\providecommand{\tp}{\mathop{\rm top}\nolimits}%
\newcommand{\Ext}{\operatorname{Ext}}
\newcommand{\pd}{\operatorname{pd}}
\newcommand{\soc}{\operatorname{soc}}
\newcommand{\op}{\operatorname{op}}

\begin{document}

\title[COMPONENTS DETERMINED BY THEIR COMPOSITION FACTORS]{AUSLANDER-REITEN COMPONENTS DETERMINED BY THEIR COMPOSITION FACTORS}

\author{ALICJA JAWORSKA}
\address{Faculty of Mathematics and Computer Science, Nicolaus Copernicus University, Chopina 12/18, 87-100 Toru\'n, Poland}
\email{jaworska@mat.uni.torun.pl}
\thanks{The research supported by the Research Grant N N201 269135 of the Polish Ministry of Science and Higher Education.}
\author{PIOTR MALICKI}
\address{Faculty of Mathematics and Computer Science, Nicolaus Copernicus University, Chopina 12/18, 87-100 Toru\'n, Poland}
\email{pmalicki@mat.uni.torun.pl}

\author{ANDRZEJ SKOWRO\'NSKI}
\address{Faculty of Mathematics and Computer Science, Nicolaus Copernicus University, Chopina 12/18, 87-100 Toru\'n, Poland}
\email{skowron@mat.uni.torun.pl}



\begin{abstract}
We provide sufficient conditions for a component of the Auslander-Reiten quiver of an artin algebra to be determined by the composition factors
of its indecomposable modules.
\end{abstract}

\maketitle

\section{Introduction and main results}
Let $A$ be an artin algebra over a commutative artin ring $R$. We denote by $\mod A$ the category of finitely generated right $A$-modules,
by $K_0(A)$ the Grothendieck group of $A$, and by $[M]$ the image of a module $M$ from $\mod A$ in $K_0(A)$. Thus, for modules $M$ and $N$
in $\mod A$, $[M]=[N]$ if and only if $M$ and $N$ have the same composition factors including the multiplicities. An interesting open problem
is to find handy criteria for two indecomposable modules $M$ and $N$ in $\mod A$ with the same composition factors to be isomorphic. It was shown
in \cite{RSS} that it is the case when $M$ does not lie on a short cycle $M\to X\to M$ of non-zero non-isomorphisms in $\mod A$ with $X$ an
indecomposable module, generalizing earlier results about directing modules proved in \cite{H}, \cite{HRi}.
In fact, it follows from \cite{HL} and \cite{RSS} that an indecomposable module $M$ in $\mod A$ lies on a short cycle $M\to X\to M$ in $\mod A$
if and only if $M$ is the middle term of a chain $Y\to M\to D\Tr Y$ of non-zero homomorphisms in $\mod A$ with $Y$ a non-projective indecomposable module.
Hence the above result from \cite{RSS} gives in fact another interpretation of a result from \cite{AR}.
An important combinatorial and homological invariant of the module category $\mod A$ of an artin algebra $A$ is its Auslander-Reiten quiver $\Gamma_A$
\cite{ARS}.
Sometimes, we may recover the algebra $A$ and the category $\mod A$ from the shape of components $\mathscr{C}$ of $\Gamma_A$ and their behaviour
in the category $\mod A$. By a component of $\Gamma_A$ we mean a connected component of the translation quiver $\Gamma_A$.

In this article we are concerned with the problem of finding handy criteria for a component $\mathscr{C}$ of the Auslander-Reiten quiver $\Gamma_A$ of
an artin algebra $A$ to be uniquely determined in $\Gamma_A$ by the composition factors of its indecomposable modules. We say that two
components $\mathscr{C}$ and $\mathscr{D}$ of $\Gamma_A$ have the same composition factors if, for any element $\mathbf{x}\in K_0(A)$, $\mathbf{x}=[M]$
for an indecomposable module $M$ in $\mathscr{C}$ if and only if $\mathbf{x}=[N]$ for an indecomposable module $N$ in $\mathscr{D}$.

In order to state the main results, we recall some concepts. For an artin algebra $A$, we denote by $\rad_A$ the Jacobson radical of $\mod A$,
generated by all non-isomorphisms between indecomposable modules in $\mod A$, and by $\rad_A^{\infty}$ the infinite radical of $\mod A$, which is the
intersection of all powers $\rad_A^i$, $i\geq 1$, of $\rad_A$. Recall that, by a result of M. Auslander \cite{A}, $\rad_A^{\infty}=0$ if and only if
$A$ is of finite representation type, that is, there are in $\mod A$ only finitely many indecomposable modules up to isomorphism. Following \cite{S3},
a component
quiver $\Sigma_A$ of $A$ is the quiver whose vertices are the components $\mathscr{C}$ of $\Gamma_A$, and two components
$\mathscr{C}$ and $\mathscr{D}$ of $\Gamma_A$ are linked in $\Sigma_A$ by an arrow $\mathscr{C}\to\mathscr{D}$ provided $\rad_A^{\infty}(X,Y)\neq 0$
for some modules $X\in\mathscr{C}$ and $Y\in\mathscr{D}$. We note that a component $\mathscr{C}$ of $\Gamma_A$ is generalized standard
in the sense of \cite{S1} if and only if $\Sigma_A$ has no loop at $\mathscr{C}$. By a short cycle in $\Sigma_A$ we mean a cycle
$\mathscr{C}\to\mathscr{D}\to\mathscr{C}$, where possibly $\mathscr{C}=\mathscr{D}$. We also mention that a component
$\mathscr{C}$ of $\Gamma_A$ lies on a short cycle $\mathscr{C}\to\mathscr{D}\to\mathscr{C}$ in $\Sigma_A$ with $\mathscr{C}\neq\mathscr{D}$
if and only if $\mathscr{C}$ has an external short path $X\to Y\to Z$ with $X$ and $Z$ in $\mathscr{C}$ and $Y$ in $\mathscr{D}$ \cite{RS}.
Recall also that a translation quiver of the form $\mathbb{Z}\mathbb{A}_{\infty}/(\tau^r)$, $r\geq 1$, is called a stable tube of rank $r$.
We note that every regular component (without projective modules and injective modules) of the Auslander-Reiten quiver $\Gamma_A$ of an artin algebra $A$ is either a stable tube or is acyclic
(without oriented cycles) of the form $\mathbb{Z}\Delta$ for an acyclic locally finite connected valued quiver $\Delta$ (see \cite{Li}, \cite{Z}).

The following theorem is the first main result of this article.
\begin{theorem} \label{thm}
Let $A$ be an artin algebra, and $\mathscr{C}$ and $\mathscr{D}$ two components of $\Gamma_A$ with the same composition factors.
Assume that $\mathscr{C}$ is not a stable tube of rank one and does not lie on a short cycle in $\Sigma_A$. Then $\mathscr{C}=\mathscr{D}$.
\end{theorem}
Therefore, the above theorem says that a generalized standard Auslander-Reiten component $\mathscr{C}$ of an artin algebra $A$ without external
short paths, different from a stable tube of rank one, is uniquely determined in $\Gamma_A$ by the composition factors of its indecomposable modules.
We point out that the assumption on $\mathscr{C}$ not being a stable tube of rank one is essential for the validity of the above theorem.
For example, if $H$ is the path algebra $K\Delta$ of a Euclidean quiver $\Delta$ over an algebraically closed field $K$, then the component quiver
$\Sigma_H$ of $H$ is acyclic and the Auslander-Reiten quiver $\Gamma_H$ of $H$ contains infinitely many pairwise
different stable tubes of rank one having the same composition factors (see \cite{R1}, \cite{SS1}).

The second main result of the article clarifies the situation in general.
\begin{theorem} \label{thm2}
Let $A$ be an artin algebra, $\mathscr{C}$ a stable tube of rank one in $\Gamma_A$ which does not lie on a short cycle in $\Sigma_A$, and $\mathscr{D}$
a component of $\Gamma_A$ different from $\mathscr{C}$ and having the same composition factors as $\mathscr{C}$.
Then there is a quotient algebra $B$ of $A$ such that the following statements hold:
\begin{itemize}
\item[(a)] $B$ is  a concealed canonical algebra.
\item[(b)] $\mathscr{C}$ and  $\mathscr{D}$ are stable tubes of a separating family of stable tubes of $\Gamma_B$.
\item[(c)] $\mathscr{D}$ is a stable tube of rank one.
\end{itemize}
\end{theorem}

Recall that a concealed canonical algebra is an algebra of the form $B=\End_{\Lambda}(T)$, where $\Lambda$ is a canonical algebra in the sense
of C. M. Ringel \cite{R2} (see also \cite{R1}) and $T$ is a multiplicity-free tilting module in the additive category $\add(\mathcal{P}^{\Lambda})$,
for the canonical decomposition $\Gamma_{\Lambda}= \mathcal{P}^{\Lambda} \vee \mathcal{T}^{\Lambda} \vee \mathcal{Q}^{\Lambda}$ of $\Gamma_{\Lambda}$,
with $\mathcal{T}^{\Lambda}$ the canonical infinite separating family of stable tubes of $\Gamma_{\Lambda}$. Then $\Gamma_B$ admits a decomposition
$\Gamma_B = \mathcal{P}^B \vee \mathcal{T}^B \vee \mathcal{Q}^B$, where the image $\mathcal{T}^B=\Hom_{\Lambda}(T, \mathcal{T}^{\Lambda})$
of the family $\mathcal{T}^{\Lambda}$ via the functor $\Hom_{\Lambda}(T, -): \mod \Lambda \rightarrow \mod B$ is an infinite separating
family of stable tubes  of $\Gamma_B$. Moreover, all but finitely many stable tubes of $\mathcal{T}^B$ have rank one and the same composition factors.
We also mention that, by a result of H. Lenzing and J. A. de la Pe\~na \cite{LP}, the class of concealed canonical algebras  coincides with the class
of artin algebras whose Auslander-Reiten quiver admits a separating family of stable tubes.

We exhibit in Section \ref{s3} examples of generalized standard stable tubes of arbitrary large rank which are not uniquely determined by the composition
factors. It would be interesting to clarify if an acyclic generalized standard regular component of the Auslander-Reiten quiver
of an artin algebra is uniquely determined by its composition factors (see Section 3 for related comments).

For basic background on the representation theory applied here we refer to \cite{ASS}, \cite{ARS}, \cite{R1}, \cite{SS1}, \cite{SS2}.

\section{Proofs of Theorems \ref{thm} and \ref{thm2}}
Let $A$ be an artin algebra over a commutative artin ring $R$. We denote by $\tau_A$ and $\tau^-_A$ the Auslander-Reiten translations $D\Tr$ and $\Tr D$,
respectively. For a module $V$ in $\mod R$, we denote by $|V|$ its length over $R$. In the proofs a crucial role will be played by the following
formulas from \cite[Proposition 4.1]{S2}, being consequences of \cite[(1.4)]{AR} (see also \cite[Corollary IV.4.3]{ARS}).

\textit{For indecomposable modules $M, N$ and $X$ in $\mod A$ with $[M]=[N]$ the following equalities hold:}

\begin{itemize}
\item[(i)] $|\!\Hom_A(X,M)| - |\!\Hom_A(M,\tau_A X)| = |\!\Hom_A(X,N)| - |\!\Hom_A(N,\tau_A X)|$,
\item[(ii)] $|\!\Hom_A(M,X)| - |\!\Hom_A(\tau^-_A X,M)| = |\!\Hom_A(N,X)| - |\!\Hom_A(\tau^-_A X,N)|$.
\end{itemize}

\noindent Let $\mathscr{C}$ and $\mathscr{D}$ be components of $\Gamma_A$ with the same composition factors and $\mathscr{C}$ does not lie
on a short cycle in $\Sigma_A$. We assume that $\mathscr{C}\neq\mathscr{D}$ and show in several steps that $\mathscr{C}$ and $\mathscr{D}$
are stable tubes of rank one of a separating family of stable tubes in the Auslander-Reiten quiver $\Gamma_B$ of a concealed canonical algebra $B$.\\

\noindent (1) $\mathscr{C}$ is a semi-regular component of $\Gamma_A$ ($\mathscr{C}$ does not contain both a projective module and an injective module).
Assume $\mathscr{C}$ contains a projective module $P$ and an injective module $I$. Since $\mathscr{C}$ and $\mathscr{D}$ have the same composition
factors, there exist modules $M$ and $N$ in $\mathscr{D}$ such that $[P]=[M]$ and $[I]=[N]$. Then we have $\Hom_A(P,M)\neq 0$ and $\Hom_A(N,I)\neq 0$,
because the top of $P$ is a composition factor of $M$, and the socle of $I$ is a composition factor of $N$.
Hence, we have in $\Sigma_A$ the short cycle $\mathscr{C}\to\mathscr{D}\to\mathscr{C}$, because $\Hom_{A}(P,M)=\rad_A^{\infty}(P,M)$ and
$\Hom_A(N,I)=\rad^{\infty}_A(N,I)$, a contradiction. Therefore, $\mathscr{C}$ is a semi-regular component of $\Gamma_A$. \\

\noindent (2) $\mathscr{C}$ is a cyclic component of $\Gamma_A$ (every module in $\mathscr{C}$ lies on an oriented cycle in $\mathscr{C}$).
Take a module $X$ in $\mathscr{C}$. It follows from our assumption that $[X]=[Y]$ for some module $Y$ in $\mathscr{D}$, and so
$X$ is not uniquely determined by $[X]$, because $\mathscr{C}\neq\mathscr{D}$. Applying \cite[Corollary 2.2]{RSS}, we conclude that we have in $\mod A$
a short cycle $X\to Z\to X$. Observe that then $Z$ belongs to $\mathscr{C}$, because $\mathscr{C}$ does not lie on a short cycle in $\Sigma_A$.
Moreover, since there is no loop at $\mathscr{C}$ in $\Sigma_A$,
 $\mathscr{C}$ is a generalized standard component of $\Gamma_A$, and hence $\rad^{\infty}_A(X,Z)=0$ and $\rad^{\infty}_A(Z,X)=0$.
Then $\Hom_A(X,Z)\neq 0$ and $\Hom_A(Z,X)\neq 0$ imply that there exist paths of irreducible homomorphisms in $\mod A$ from $X$ to $Z$ and from $Z$ to $X$
(see \cite[Proposition V.7.5]{ARS}), and consequently an oriented cycle in $\mathscr{C}$ passing through $X$ and $Z$. Hence, $\mathscr{C}$ is a cyclic
component. \\

\noindent (3) $\mathscr{C}$ is a ray tube (obtained from a stable tube by a finite number (possibly empty) of ray insertions) or
a coray tube (obtained from a stable tube by a finite number (possibly empty) of coray insertions) in the sense of \cite[(4.5)]{R1}
(see also \cite[XV.2]{SS2}). This is a direct consequence of \cite[(2.6)]{Liu}, since by (1) and (2) $\mathscr{C}$ is semi-regular with oriented
cycles. \\

\noindent (4) We may assume (without loss of generality) that $\mathscr{C}$ is a ray tube, hence without injective modules. Let
$\ann_A(\mathscr{C})$ be the annihilator of $\mathscr{C}$ in $A$, that is, the intersection of the annihilators $\ann_A(X)=\{a\in
A\mid Xa=0\}$ of all modules $X$ in $\mathscr{C}$, and $B=A/\ann_A(\mathscr{C})$. Then $\mathscr{C}$ is a faithful
component of $\Gamma_B$. Since $\mathscr{C}$ does not lie on a short cycle in $\Sigma_A$, we conclude that $\mathscr{C}$ is
without external short paths \cite{RS}, that is, there are no paths $U\to V\to W$ in $\mod A$ with $U$ and $W$ in $\mathscr{C}$
but $V$ not in $\mathscr{C}$. Then it follows from \cite[Theorem 2]{JMS} that $B$ is an almost concealed canonical algebra and
$\mathscr{C}$ is a faithful ray tube of a separating family $\mathcal{T}^{B}$ of ray tubes of $\Gamma_B$. Recall that then
there exists a canonical algebra $\Lambda$ (in the sense of C. M. Ringel \cite{R1}, \cite{R2}) such that $B=\End_{\Lambda}(T)$ for a
tilting module $T$ in the additive category $\add(\mathcal{P}^{\Lambda}\cup \mathcal{T}^{\Lambda})$ of
$\mathcal{P}^{\Lambda}\cup \mathcal{T}^{\Lambda}$, for the canonical decomposition $\Gamma_{\Lambda}=\mathcal{P}^{\Lambda}
\vee \mathcal{T}^{\Lambda} \vee \mathcal{Q}^{\Lambda}$ of $\Gamma_{\Lambda}$ with $\mathcal{T}^{\Lambda}$ the canonical
separating family of stable tubes. By general theory (see \cite{LP}, \cite{LS}, \cite{R1}, \cite{R2}, \cite{S3a}),
$\Gamma_B$ admits a decomposition
\[\Gamma_B=\mathcal{P}^B \vee \mathcal{T}^B \vee \mathcal{Q}^B,\]
where $\mathcal{T}^B$ is a family of ray tubes separating $\mathcal{P}^B$ from $\mathcal{Q}^B$ (in the sense of \cite{R2}). In particular,
$\mathcal{T}^B$ is an infinite family of pairwise orthogonal generalized standard ray tubes, $\Hom_B(\mathcal{T}^B, \mathcal{P}^B)$ $=0$,
$\Hom_{B}(\mathcal{Q}^B, \mathcal{T}^B)=0$,
and $\Hom_B(\mathcal{Q}^B, \mathcal{P}^B)=0$. In fact, since $\mathscr{C}$ is a faithful ray tube of $\mathcal{T}^B$, all ray tubes of $\mathcal{T}^{B}$
except $\mathscr{C}$ are  stable tubes. Moreover, the separation property of $\mathcal{T}^B$ implies that $\Hom_{B}(\mathscr{X}, \mathscr{C})\neq 0$
for any component $\mathscr{X}$ from $\mathcal{P}^B$ and $\Hom_{B}(\mathscr{C}, \mathscr{Y}) \neq 0$ for any component $\mathscr{Y}$ from $\mathcal{Q}^B$.
Moreover, we note that $\mathcal{Q}^B$ contains all indecomposable injective $B$-modules. \\

\noindent (5) $\mathscr{D}$ is a component of $\Gamma_B$. Write $A=P'\oplus P''$ where the simple summands of $P'/\rad P'$ are exactly the simple
composition factors of modules in $\mathscr{C}$. Denote by $t_{P''}(A)$ the ideal of $A$ generated by the images of all homomorphisms
in $\mod A$ from $P''$ to $A$.
Since $\mathscr{C}$ is a semi-regular component of $\Gamma_A$ without external short paths,
it follows from arguments in  \cite[Section 1]{RS} that $\End_A(P')\cong A/t_{P''}(A)$ and $t_{P''}(A)=\ann_A(\mathscr{C})$. Observe that $1_A=e+f$ for
orthogonal idempotents $e$ and $f$ in $A$ with $P'=eA$ and $P''=fA$, and consequently $\End_A(P')\cong eAe$ and $t_{P''}(A)=AfA$. Clearly, then
$B=A/\ann_A(\mathscr{C})\cong eAe$.
On the other hand, since $\mathscr{D}$ has the same composition factors as $\mathscr{C}$, we have $Nf=\Hom_A(fA,N)=\Hom_A(P'',N)=0$, and consequently
$N\ann_A(\mathscr{C})=N(AfA)=(Nf)A=0$, for any module $N$ in $\mathscr{D}$.
This shows that $\mathscr{D}$ is a component of $\Gamma_B$. \\

\noindent (6) $\mathscr{D}$ is a component of $\mathcal{T}^B$. Assume $\mathscr {D} \notin \mathcal{T}^B$. Fix a stable tube
$\mathcal{T}^{\ast}$ of $\mathcal{T}^B$ of rank one, which is different from $\mathscr{C}$. By general theory (\cite{LP},
\cite{LS}, \cite{R2}) B is a tubular (branch) extension of a concealed canonical algebra $C$ such that $\Gamma_C=\mathcal{P}^C
\vee \mathcal{T}^{C} \vee \mathcal{Q}^{C}$, where $\mathcal{T}^C$ is a separating family of stable tubes,
$\mathcal{P}^B=\mathcal{P}^C$, $\mathscr{C}$ is obtained from a stable tube $\mathcal{T}$ of $\mathcal{T}^C$ by a finite number
(possibly empty) of ray insertions and the remaining tubes of $\mathcal{T}^C$ and $\mathcal{T}^B$ coincide
($\mathcal{T}^C\backslash \mathcal{T}=\mathcal{T}^B \backslash \mathscr{C}$). Clearly, $C$ is a quotient algebra of $B$.

Let $M$ be a module in $\mathscr{C}$ which lies in $\mathcal{T}$. In particular, the composition factors of $M$ are $C$-modules.
Take a module $N \in \mathscr{D}$ such
that $[M]=[N]$. Assume $\mathscr{D} \in \mathcal{Q}^B$. Since $[M]=[N]$ there exists a projective module $P \in
\mathcal{P}^B=\mathcal{P}^C$ such that $\Hom_B(P,N)\neq 0$. By the separation property of $\mathcal{T}^B$ we have $\Hom_B(X,N)\neq 0$
for some module $X \in \mathcal{T}^{\ast}$. Then, applying the formula (i), we obtain
\[0=|\!\Hom_A(X,M)| - |\!\Hom_A(M,\tau_A X)| = |\!\Hom_A(X,N)| - |\!\Hom_A(N,\tau_A X)|\] \[ =|\!\Hom_A(X,N)|>0,\]
since $M$ and $X$ belong to orthogonal tubes of $\mathcal{T}^B$ and $\Hom_B(\mathcal{Q}^B, \mathcal{T}^B)=0$. Dually, if $\mathcal{D}
\in \mathcal{P}^B$, then there exists an injective module $I$ in $\mathcal{Q}^B$ such that $\Hom_B(N,I)\neq 0$. By the separation
property of $\mathcal{T}^B$, we have $\Hom_B(N,Y) \neq 0$ for some module $Y \in \mathcal{T}^{\ast}$. Then, by the formula (ii), we have
\[0= |\!\Hom_A(M,Y)| - |\!\Hom_A(\tau^-_A Y,M)| = |\!\Hom_A(N,Y)| - |\!\Hom_A(\tau^-_A Y,N)|\] \[=|\!\Hom_A(N,Y)|>0,\]
since $M$ and $Y$ belong to orthogonal tubes of $\mathcal{T}^B$ and $\Hom_B(\mathcal{T}^B, \mathcal{P}^B)=0$. The above contradictions
show that $\mathscr{D} \in \mathcal{T}^B$.\\

\noindent (7) $\mathcal{T}^B$ is a family of stable tubes. Assume $\mathscr{C}$ contains a projective module $P$. Take an
indecomposable module $Y$ in $\mathscr{D}$ with $[P]=[Y]$. Then the top of $P$ is a composition factor of $Y$ and hence $\Hom_B(P,Y)\neq 0$.
Therefore, $\Hom_B(\mathscr{C},\mathscr{D})\neq 0$ which contradicts the fact
that $\mathscr{C}$ and $\mathscr{D}$ are orthogonal. We conclude that $\mathscr{C}$ is a stable tube of $\mathcal{T}^B$. Clearly,
then $\mathcal{T}^B$ is a separating family of stable tubes of $\Gamma_B$, and consequently $B$ is a concealed canonical algebra,
by \cite{LP}. \\

\noindent (8) $\mathscr{C}$ and $\mathscr{D}$ are stable tubes of rank one. Since $\mathscr{C}$ and $\mathscr{D}$ belong to the
separating family $\mathcal{T}^B$ of stable tubes of $\Gamma_B$, we know that $\mathscr{C}$ and $\mathscr{D}$ are orthogonal,
generalized standard, and without external short paths. In particular, $\mathscr{C}$ and $\mathscr{D}$ do not lie on short
cycles in $\Sigma_B$. Then, applying \cite[Lemmas 3.1 and 3.3]{S2}, we conclude that $\mathscr{C}$ and $\mathscr{D}$ consist of modules
which do not lie on infinite short cycles in $\mod B$. Assume $\mathscr{C}$ is of rank $r\geq 2$. Take a module $X$ lying on the
mouth of $\mathscr{C}$ ($X$ has one immediate predecessor and one immediate successor in $\mathscr{C}$). Then, by \cite[Corollary
4.4]{S2}, $X$ is uniquely determined by $[X]$, which contradicts the fact that $[X]=[Y]$ for some module $Y$ in $\mathscr{D}$ and
$\mathscr{C}\neq\mathscr{D}$. Therefore, $\mathscr{C}$ is of rank one. Applying the same arguments, we conclude that $\mathscr{D}$ is
also of rank one.

Summing up, the proofs of Theorems \ref{thm} and \ref{thm2} are provided.

\section{Examples} \label{s3}

Let $K$ be an algebraically closed field and $Q$ be a finite quiver. For any arrow $\alpha \in Q$, by $s(\alpha)$ and $t(\alpha)$
we mean the source and the target of $\alpha$, respectively. By $KQ$ we denote the path algebra of $Q$.
Recall that, if the quiver $Q$ is acyclic, then $KQ$ is a hereditary algebra \cite{ASS}. For a finite dimensional algebra
$H$ over $K$, we denote by $T(H)$ the trivial extension algebra of  $H$ by its duality $H$-$H$-bimodule $D(H)=\Hom_K(H,K)$.
Recall that $T(H)=H \oplus D(H)$ as $K$-vector space and the multiplication in $T(H)$ is given by $(a,f)(b,g)=(ab, ag+fb)$
for $a,b \in H$ and $f,g \in D(H)$. Then $T(H)$ is a symmetric algebra and $H$ is the quotient algebra of $T(H)$ by the ideal $D(H)$.

For a natural number $n\geq 4$, $Q_{n}$ will be the quiver of the following form:
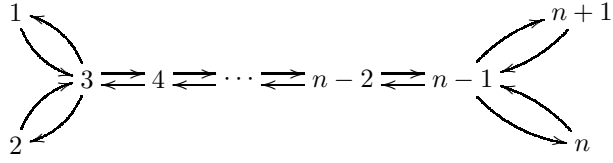
\[\xymatrix@C=16pt@R=12pt{1 \ar@/_0.5pc/[dr] & &&&& &
  n+1 \ar@/^0.5pc/[ld]\\
           & 3 \ar@<0.5ex>[r] \ar@/_0.5pc/[lu] \ar@/^0.5pc/[ld] & 4 \ar@<0.5ex>[l] \ar@<0.5ex>[r]
           &  \ar@<0.5ex>[l] \cdots \ar@<0.5ex>[r]&  \ar@<0.5ex>[l]  n-2  \ar@<0.5ex>[r]& \ar@<0.5ex>[l] n-1  \ar@/^0.5pc/[ru]  \ar@/_0.5pc/[rd] \\
           2 \ar@/^0.5pc/[ur] & &&& && n \ar@/_0.5pc/[lu] }\]
Each arrow in $Q_{n}$ will be named either by $\alpha$ or by $\beta$ in such a way that an arrow which starts in the vertex $3$ and terminates in
the vertex $1$ is $\alpha$, and $s(\alpha)=t(\beta)$, $t(\alpha)=s(\beta)$, for all arrows $\alpha$ and $\beta$.
Let $I_{n}$ be the admissible ideal in the path algebra $KQ_{n}$ generated by all paths $\alpha\beta,\beta\alpha$ such that $s(\alpha\beta)\neq t(\alpha\beta),
s(\beta\alpha)\neq t(\beta\alpha)$, and all commutativity relations $\omega_{1}-\omega_{2}$, where $\omega_{1},\omega_{2}$ are all paths of length $2$ in $Q_n$ such that
their source and target coincide with the vertex $i$, for all $i \in \{3,\ldots,n-1\}$. Then by $\Lambda _n$ we denote the quotient algebra $KQ_n/ I_n$.\\

We consider now the quiver $\Delta_{n}$ of Euclidean type $\widetilde{\mathbb{D}}_{n}$, for any $n\geq 4$, defined in the following way.
If $n$ is an odd number, then $\Delta_{n}$ is of the form:
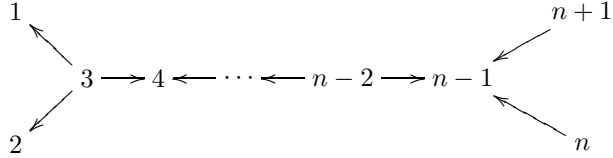
\[ \xymatrix@C=16pt@R=12pt{1  & &&&&& {n+1} \ar[ld]\\
           & 3 \ar[r] \ar[lu] \ar[ld] & 4
           &  \ar[l]  \cdots &\ar[l] n-2 \ar[r]&  {n-1}  \\
           2  & &&&&& {n} \ar[lu]}\]
 and similarly, for an even number $n$,  the quiver $\Delta_{n}$ is of the form:
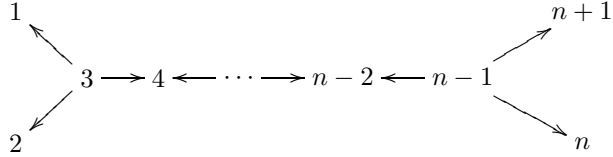
\[\xymatrix@C=16pt@R=12pt{1  & &&&&& {n+1} \\
           & 3 \ar[r] \ar[lu] \ar[ld] & 4
           &  \ar[l] \cdots \ar[r]&  n-2  & \ar[l] {n-1}  \ar[ru]  \ar[rd] \\
           2  & &&&& &{n} }\]
(in particular, all maximal subquivers of type $\mathbb{A}_{n-1}$ of $Q_n$ have alternate orientation of arrows).

Let $H_n$ be the path algebra $K\Delta_n$ and $H_n^{\ast}$ the path
algebra $K\Delta^{\ast}_n$, where $\Delta^{\ast}_n$ is the opposite
quiver of $ \Delta_n$.
Note that $\Delta_{n}$ is a subquiver of $Q_{n}$ given by the arrows $\alpha$ and $\Delta_{n}^{*}$ is a subquiver of $Q_{n}$ given by the arrows $\beta$.
Moreover, observe that $\Lambda_{n}$ is the trivial extension algebra $T(H_{n})$ of  $H_{n}$ and the trivial extension algebra
$T(H^{*}_{n})$ of $H^{*}_{n}$. In particular, $H_n$ and $H_n^{\ast}$ are quotient algebras of  $\Lambda_n$.  \\

Assume now that $n \geq 4$ is an odd number. For each arrow $\alpha$ in $\Delta_n$ such that $s(\alpha)=i$ and $t(\alpha)\in \{i-1,i+1\}$, for some
$i \in \{3,\ldots,n\}$, we put $\alpha_l$ instead of $\alpha$, where $l$ is given by the formula:
\[l=\begin{cases}
\frac{i-1}{2}; \qquad \quad \alpha: i \rightarrow i+1\\
n-\frac{i+1}{2};\quad \alpha: i \rightarrow i-1.
\end{cases}
\]
Observe that $l \in \{1,\ldots,n-2\}$. We define the family of indecomposable representations $F_{\alpha_1}, \ldots,F_{\alpha_{n-2}}$  of $H_n$ over $K$:
\begin{enumerate}
\item[$\bullet$] $F_{\alpha_l}$ for $l \notin \{\frac{n-1}{2}, n-2\}:$
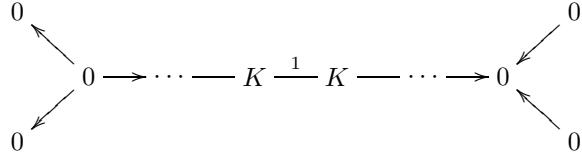
\[\xymatrix@C=16pt@R=12pt{0 & &&&&&& 0 \ar[ld]\\
           & 0 \ar[r] \ar[lu] \ar[ld]
           & \cdots \ar@{-}[r]& K \ar@{-}[r]^(0.5)1 & K \ar@{-}[r]& \cdots \ar[r]& 0  \\
           0  & &&&&&& 0 \ar[lu]}\]
           where $K$ stands  in  the vertices $s(\alpha_l), t(\alpha_l)$, zero space elsewhere\\
           (here by $\xymatrix@C=16pt@R=12pt{\ar@{-}[r]&}$ we mean $\xymatrix@C=16pt@R=12pt{\ar[r]&}$ or $\xymatrix@C=16pt@R=12pt{&\ar[l]}$);
\item[$\bullet$] $F_{\alpha_{\frac{n-1}{2}}}:$
 
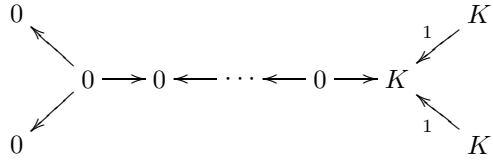
\[ \xymatrix@C=16pt@R=12pt{0 & &&&&& K \ar[ld]_(0.5)1\\
           & 0 \ar[r] \ar[lu] \ar[ld]
           &  0  & \ar[l]  \cdots & \ar[l] 0 \ar[r]& K \\
           0  & &&&&& K \ar[lu]^(0.5)1}\]
            where $K$ stands in the vertices  $n-1,n,n+1$, zero space elsewhere;
\item[$\bullet$] $F_{\alpha_{n-2}}:$
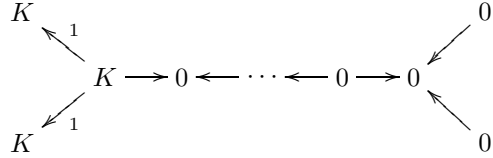
\[  \xymatrix@C=16pt@R=12pt{K & &&&&& 0 \ar[ld]\\
           & K \ar[r] \ar[lu]_(0.5)1 \ar[ld]^(0.5)1
           & 0 & \ar[l] \cdots  & 0 \ar[l] \ar[r]& 0 \\
           K  & &&&&& 0 \ar[lu]}\]
            where $K$ stands in the vertices $1,2,3$, zero space elsewhere.
\end{enumerate}

Let $E_l=F_{\alpha_l}$ for $l \in \{1,\ldots,n-2\}$. Obviously $E_1$, \ldots, $E_{n-2}$ are pairwise orthogonal bricks. Direct calculation
shows that $\tau_{H_n}E_{l+1}=E_l$ if $l \in \{1,\ldots,n-3\}$ and $\tau_{H_n}E_1=E_{n-2}$. Moreover, $\Ext ^2_{H_n}(E_{r},E_p)=0$
for any $r,p \in \{1,\ldots,n-2\}$, because $H_n$ is a hereditary algebra. It allows us to state that $E_1$, \ldots, $E_{n-2}$ form the
mouth of a standard stable tube $\mathcal{T}$ of rank $n-2$ in $\Gamma_{H_n}$ (see \cite{R1},\cite{SS1}). Since $\pd_{H_n}X \leq 1$ for any
$H_n$-module $X$ in $\mathcal{T}$, it follows from \cite[Proposition 1.1]{S4} that $\mathcal{T}$ is also a component  of the Auslander-Reiten
quiver $\Gamma_{\Lambda_n}$.\\

Analogously, let $E^*_1, E^*_2,\ldots,E^*_{n-2}$ be the indecomposable $H^{\ast}_n$-modules, where the indices $l$ are given in such a way that,
for any $l \in \{1,\ldots,n-2\}$, $E_l$ and $E_l^*$ have the same composition factors in $\mod\Lambda_n$ including the multiplicities. It is easy to see
that these modules form the mouth of  a stable tube $\mathcal{T}^*$ of rank $n-2$ in $\Gamma_{H^*_n}$ such
that $\tau_{H^{\ast}_n}E^{\ast}_{l}=E^{\ast}_{l+1}$ for $l \in \{1,\ldots, n-3\}$ and $\tau_{H_n^{\ast}}E^{\ast}_{n-2}=E^{\ast}_{1}$.
Using once more  \cite[Proposition 1.1]{S4} we get that  $\mathcal{T}^{\ast}$ is also a component of the Auslander-Reiten quiver $\Gamma_{\Lambda_n}$.
Note that $\tp(E^*_l)=\soc(E_l)$ and $\tp(E_l)=\soc(E_l^*)$ in $\mod \Lambda_n$, for any $l \in \{1, \ldots,n-2\}$. Therefore, $\mathcal{T}$ has an external
short path $E_l \rightarrow E_l^* \rightarrow E_l$ in $\mod{\Lambda_n}$, which implies existence of a short cycle
$\mathcal{T} \rightarrow \mathcal{T}^* \rightarrow \mathcal{T}$ in $\Sigma_{\Lambda_n}$. Observe also that  $\mathcal{T}$ and $\mathcal{T}^{\ast}$
have the same composition factors since $[E_l]=[E^*_l]$ for all $l \in \{1,\ldots,n-2\}$. Moreover, $\mathcal{T}$ and $\mathcal{T}^*$ are
generalized standard stable tubes in $\Gamma_{\Lambda_n}$ since they are generalized standard in $\Gamma_{H_n}$ and $\Gamma_{H_n^{\ast}}$,
respectively (see for example \cite[Chapter X]{SS1}).

Assume $n\geq 4$ is an even number. For each arrow $\alpha$ in $H_n$ such that $s(\alpha)=i$ and $t(\alpha)\in \{i-1,i+1\}$, for some
 $i \in \{3,\ldots,n-1\}$, we define the index $l$ in the previous way. Similarly, we define the indecomposable representations
 $F_{\alpha_1},\ldots, F_{\alpha_{n-2}}$ of $H_n$ over $K$:\\

\begin{enumerate}
\item[$\bullet$] $F_{\alpha_l}$ for $l \notin \{\frac{n-2}{2}, n-2\}:$
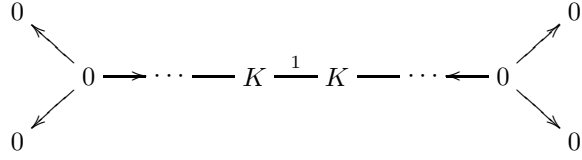
\[\xymatrix@C=16pt@R=12pt{0 & &&&&&& 0 \\
           & 0 \ar[r] \ar[lu] \ar[ld]
           & \cdots \ar@{-}[r]& K \ar@{-}[r]^(0.5)1 & K  \ar@{-}[r]&\cdots & 0 \ar[ru] \ar[rd] \ar[l] \\
           0 & &&&&&& 0}\]
           where $K$ stands in the vertices $s(\alpha_l), t(\alpha_l)$, zero space elsewhere\\
            (here by $\xymatrix@C=16pt@R=12pt{\ar@{-}[r]&}$ we mean $\xymatrix@C=16pt@R=12pt{\ar[r]&}$ or $\xymatrix@C=16pt@R=12pt{&\ar[l]}$);
\item[$\bullet$] $F_{\alpha_{\frac{n-2}{2}}}:$

 
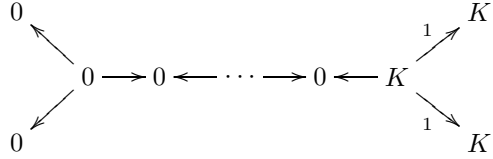
\[ \xymatrix@C=16pt@R=12pt{0 & &&&&& K \\
           & 0 \ar[r] \ar[lu] \ar[ld]
           & 0 &\ar[l] \cdots\ar[r] &  0 & K \ar[l]\ar[ru]^(0.5)1\ar[rd]_(0.5)1\\
           0  & &&&&& K }\]
            where $K$ stands in the vertices $n-1,n,n+1$, zero space elsewhere;
\item[$\bullet$] $F_{\alpha_{n-2}}:$

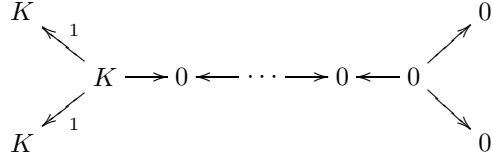
\[  \xymatrix@C=16pt@R=12pt{K & &&&&& 0 \\
           & K \ar[r] \ar[lu]_(0.5)1 \ar[ld]^(0.5)1
           & 0 &\ar[l] \cdots\ar[r] &  0   & 0  \ar[ru] \ar[rd] \ar[l] \\
           K  & &&&&& 0}\]
            where $K$ stands in the vertices $1,2,3$,  zero space elsewhere.
\end{enumerate}
As before the modules $E_l=F_{\alpha_l}$, $l \in \{1,\ldots,n-2\}$, form the mouth of a stable tube $\mathcal{T}$ of rank $n-2$ in $\Gamma_{H_n}$,
in such a way that $\tau_{H_n}E_{l+1}=E_l$  for $l \in \{1,\ldots, n-3\}$ and $\tau_{H_n}E_1=E_{n-2}$. Similarly, let $\mathcal{T}^{\ast}$ be the
stable tube  of rank $n-2$ in $\Gamma_{H_n^{\ast}}$ whose mouth consists of the modules $E^*_1, E^*_2,\ldots,E^*_{n-2}$, where the indices $l$
are given in such a way that, for any $l \in \{1,\ldots,n-2\}$, $E_l$ and $E_l^*$ have the same composition factors
and $\tp(E^*_l)=\soc(E_l)$, $\tp(E_l)=\soc(E_l^*)$ in $\mod \Lambda_n$.
Therefore, there is a short cycle $\mathcal{T} \rightarrow \mathcal{T}^* \rightarrow \mathcal{T}$
in $\Sigma_{\Lambda_n}$. Moreover, $\mathcal{T}$ and $\mathcal{T}^{\ast}$ are generalized standard components in $\Gamma_{\Lambda_n}$.\\

Summing up, we have proved that, for an arbitrary $m \geq 2$, the Auslander-Reiten quiver $\Gamma_{\Lambda_{m+2}}$
of $\Lambda_{m+2}$ contains a generalized standard stable tube of rank $m$ which
is not uniquely determined by its composition factors.\\

We end this section with comments concerning acyclic generalized standard Auslander-Reiten components. It has been proved
in \cite[Corollaries 2.4 and 3.3]{S1} that every acyclic generalized standard component $\mathscr{C}$ of the
Auslander-Reiten quiver $\Gamma_A$ of an artin algebra $A$ is of the form $\mathbb{Z}\Delta$ for a finite acyclic connected
valued quiver $\Delta$ with at least three vertices, $B=A / \ann_A(\mathscr{C})$ is a tilted algebra of the form $\End_H(T)$,
for some wild hereditary artin algebra $H$ and a regular tilting $H$-module, and $\mathscr{C}$ is the connecting component
$\mathscr{C}_T$ of $\Gamma_B$ determined by $T$. Moreover, C. M. Ringel proved in \cite{Ri} that, for any connected wild
hereditary artin algebra $H$ whose ordinary valued quiver has at least three vertices, there exists a multiplicity-free regular tilting module
$T$ in $\mod H$, and  consequently the connecting component $\mathscr{C}_T$ of the Auslander-Reiten quiver $\Gamma_B$
of the associated tilted algebra $B=\End_H(T)$ is an acyclic generalized standard faithful regular component of $\Gamma_B$.
We refer also to \cite{KS} for constructions of tilted algebras having regular connecting components with arbitrary large
composition factors.

Let $K$ be an algebraically closed field, $Q$ an arbitrary connected acyclic wild quiver with at least three vertices,
and $H=KQ$. Then it follows from \cite[Corollary 4]{KS} that, there are infinitely many pairwise non-isomorphic tilted algebras
$B=\End_H(T)$, for multiplicity-free regular tilting modules $T$ in $\mod H$, such that the connecting component
$\mathscr{C}_T$ determined by $T$ is regular and without simple modules. Take such a tilted algebra $B=\End_H(T)$ and consider
the trivial extension algebra $\Lambda=T(B)$ of $B$ by the $B$-$B$-bimodule $D(B)=\Hom_K(B,K)$. Then it follows from
\cite[Section 5]{EKS} that the Auslander-Reiten quiver $\Gamma_{\Lambda}$ of $\Lambda$ consists of two acyclic generalized standard
regular sincere components $\mathscr{C}=\mathscr{C}_T$ and $\mathscr{D}$, having sections of type $\Delta=Q^{\op}$,
and infinitely many components whose stable parts are of the form $\mathbb{Z}\mathbb{A}_{\infty}$. However, it is not clear
if $\mathscr{C}$ and $\mathscr{D}$ may have the same composition factors. It would be interesting to know if such a situation may
occur.

\bibliographystyle{amsplain}

\begin{thebibliography}{10}

\bibitem {ASS} I. Assem, D. Simson and A. Skowro\'{n}ski, \textit{Elements of the Representation Theory of Associative Algebras 1:
Techniques of Representation Theory}, London Math. Soc. Stud. Texts \textbf{65} (Cambridge Univ. Press, 2006).
%
\bibitem {A} M. Auslander, \textit{Representation theory of artin algebras II},
Comm. Algebra \textbf{1} (1974), 269--310.

\bibitem {AR} M. Auslander and I. Reiten, \textit{Modules determined by their composition factors},
Illinois J. Math. \textbf{29} (1985), 280--301.

\bibitem {ARS} M. Auslander, I. Reiten and S. O. Smal{\o}, \textit{Representation Theory of Artin Algebras}, Cambridge Stud.
Adv. Math. \textbf{36} (Cambridge Univ. Press, 1995).

\bibitem {EKS} K. Erdmann, O. Kerner and A. Skowro\'nski, \textit{Self-injective algebras of wild tilted type},
J. Pure Appl. Algebra \textbf{149} (2000), 127--176.

\bibitem {H} D. Happel, \textit{Composition factors of indecomposable modules},
Proc. Amer. Math. Soc. \textbf{86} (1982), 29--31.

\bibitem {HL} D. Happel and S. Liu, \textit{Module categories without short cycles are of finite type},
Proc. Amer. Math. Soc. \textbf{120} (1994), 371--375.

\bibitem {HRi} D. Happel and C. M. Ringel, \textit{Tilted algebras},
Trans. Amer. Math. Soc. \textbf{274} (1982), 399--443.

\bibitem {JMS} A. Jaworska, P. Malicki and A. Skowro\'nski, \textit{On Auslander-Reiten components of algebras without external short paths},
J. Lond. Math. Soc., in press.

\bibitem {KS} O. Kerner and A. Skowro\'nski, \textit{Quasitilted one-point extensions of wild hereditary algebras},
J. Algebra \textbf{244} (2001), 785--827.

\bibitem {LP} H. Lenzing and J. A. de la Pe\~na, \textit{Concealed-canonical algebras and separating tubular families},
Proc. Lond. Math. Soc. (3) \textbf{78} (1999), 513--540.

\bibitem {LS} H. Lenzing and A. Skowro\'nski, \textit{Quasi-tilted algebras of canonical type},
Colloq. Math. \textbf{71} (1996), 161--181.

\bibitem {Li} S. Liu, \textit{The degrees of irreducible maps and the shapes of the components of the Auslander-Reiten quivers},
J. Lond. Math. Soc. (2) \textbf{45} (1992), 32--54.

\bibitem {Liu} S. Liu, \textit{Semi-stable components of an Auslander-Reiten quiver},
J. Lond. Math. Soc. (2) \textbf{47} (1993), 405--416.

\bibitem {RS} I. Reiten and A. Skowro\'nski, \textit{Sincere stable tubes},
J. Algebra \textbf{232} (2000), 64--75.

\bibitem {RSS} I. Reiten, A. Skowro\'{n}ski and S. O. Smal{\o}, \textit{Short chains and short cycles of modules},
Proc. Amer. Math. Soc. \textbf{117} (1993), 343--354.

\bibitem {R1} C. M. Ringel, \textit{Tame Algebras  and Integral Quadratic Forms},  Lecture Notes in Math. \textbf{1099}
(Springer Verlag, 1984).

\bibitem {Ri} C. M. Ringel, \textit{The regular components of the Auslander-Reiten quiver of a tilted algebra},
Chin. Ann. Math. Ser. B \textbf{9} (1988), 1--18.

\bibitem {R2} C. M. Ringel, \textit{The canonical algebras}, with an appendix by W. Crawley-Boevey, Topics in Algebra,
Banach Center Publ. \textbf{26},  Part 1, (PWN Warsaw, 1990), 407--432.

\bibitem {SS1} D. Simson and  A. Skowro\'{n}ski, \textit{Elements of the Representation Theory of Associative Algebras 2:
Tubes and Concealed Algebras of Euclidean Type}, London Math. Soc. Stud. Texts \textbf{71} (Cambridge Univ. Press, 2007).

\bibitem {SS2} D. Simson and  A. Skowro\'{n}ski,  \textit{Elements of the Representation Theory of Associative Algebras 3:
Representation-Infinite Tilted Algebras}, London Math. Soc. Stud. Texts \textbf{72} (Cambridge University Press, 2007).

\bibitem {S1} A. Skowro\'nski, \textit{Generalized standard Auslander-Reiten components},
J. Math. Soc. Japan \textbf{46} (1994), 517--543.

\bibitem {S2} A. Skowro\'nski, \textit{On the composition factors of periodic modules},
J. Lond. Math. Soc. (2) \textbf{49} (1994), 477--492.

\bibitem {S3} A. Skowro\'nski, \textit{Cycles in module categories},
Finite Dimensional Algebras and Related Topics, NATO ASI Series, Series C: Math. and Phys. Sciences \textbf{424}
(Kluwer Acad. Publ., 1994), 309--345.

\bibitem{S3a} A. Skowro\'nski, \textit{On omnipresent tubular families of
modules},  Representation Theory of Algebras, Canad. Math. Soc.
Conf. Proc. 18, Amer. Math Soc., Providence, RI, 1996, 641--657.

\bibitem {S4} A. Skowro\'nski, \textit{A construction of complex syzygy periodic modules over symmetric algebras},
Colloq. Math. \textbf{103} (2005), 61--69.

\bibitem {Z} Y. Zhang, \textit{The structure of stable components}, Canad. J. Math. \textbf{43} (1991), 652--672.

\end{thebibliography}

\end{document}